\documentclass[mathpazo]{cicp}
\usepackage{graphicx} % Required for inserting images
\usepackage{amssymb}
\usepackage{amsmath}
\usepackage{amsfonts}

\title{ARDO: A Weak Formulation Deep Neural Network Method for Elliptic and Parabolic PDEs Based on Random Differences of Test Functions}
\author{Andrew Qing He\affil{1}, Wei Cai\affil{1}}
\address{\affilnum{1}\ Dept. of Mathematics, Southern Methodist University, Dallas, TX 75275.}
\date{September 2025}

\begin{document}
\begin{abstract}
We propose ARDO method for solving PDEs and PDE-related problems with deep learning techniques. This method uses a weak adversarial formulation but transfers the random difference operator onto the test function. The main advantage of this framework is that it is fully derivative-free with respect to the solution neural network. This framework is particularly suitable for Fokker-Planck type second-order elliptic and parabolic PDEs.
\end{abstract}

\maketitle

\section{Introduction}
We propose \textbf{adjoint-random-difference operator (ARDO)}\footnote{We intentionally use dashes to connect the first three words to emphasize that ARDO is not a random difference operator (RDO) in the usual sense -- ARDO never computes any kind of differences or differentials on the solution network.} method, a sibling method of deep random difference method (DRDM) proposed in \cite{drdm}. Similar to SOC-MartNet \cite{socmartnet, dfsocmartnet} and DRDM method, our approach leverages a random difference method framework with a weak adversarial training \cite{wan, weran} framework but introduces a key innovation: \textit{transferring the random difference operation to the test function}. This construction yields a significant advantage by being fully derivative-free with respect to the solution neural network, thereby enhancing computational efficiency and stability during training.

ARDO is particularly suitable for solving second-order elliptic and parabolic Fokker-Planck type PDEs. By eliminating the need for higher-order derivatives of the solution network, our method mitigates the common challenges associated with gradient computation in high dimensions, offering a robust and scalable alternative for complex stochastic problems. Besides, the weak formulation provides a natural way to impose the boundary conditions, which relieves us from the troublesome trade-off between the residual loss and the boundary value loss, i.e. adjusting hyper-parameters like $\beta$ in \cite{deepritz} and $\alpha_{\text{mart}}$, $\alpha_{\text{bdry}}$ in \cite{martnet} by hand.

\section{Method}
\subsection{Weakly Nonlinear Elliptic PDEs}
Suppose $\Omega \subseteq \mathbb{R}^{n}$ is an open region with a piecewise-smooth boundary $\partial \Omega$. We consider the following elliptic or hypoelliptic Fokker-Planck style PDE
\begin{equation} \label{pg2GsuM}
-\frac{1}{2}\sum _{i,j=1}^{n}\frac{\partial ^{2}( a_{ij} f)}{\partial x_{i} \partial x_{j}} +\sum _{j=1}^{n}\frac{\partial }{\partial x_{i}}( b_{j} f) +R( f, x, t) =0,\quad x\in \Omega 
\end{equation}
with boundary conditions 
\begin{equation} \label{BcSlwdv}
\begin{cases}
f=g & x\in \partial \Omega _{\text{D}}\\
\displaystyle \sum _{j=1}^{n} \nu _{j}\left[\frac{1}{2}\sum _{i=1}^{n}\frac{\partial }{\partial x_{i}} (a_{ij} f)-b_{j} f\right] =\phi  & x\in \partial \Omega _{\text{N}}
\end{cases}
\end{equation}
where $a_{ij}$ and $b_j$ are all functions of $x$ in $\Omega$, and $R$ is an $\mathbb{R} \times \Omega\to \mathbb{R}$ source function that may be nonlinear. $\partial \Omega _{\text{D}} ,\partial \Omega _{\text{N}}$ is a partition of $\partial \Omega $ with $\partial \Omega _{\text{D}} \neq \emptyset $ (so that the solution of the equation is guaranteed to be unique). 

Suppose there is a smooth test function $\rho $ that vanishes on $\partial \Omega _{\text{D}}$. Then 
\begin{equation} \label{R4NtoMw}
\int _{\Omega }\left[ -\frac{1}{2} \rho \sum _{i,j=1}^{n}\frac{\partial ^{2}( a_{ij} f)}{\partial x_{i} \partial x_{j}} +\rho \sum _{j=1}^{n}\frac{\partial }{\partial x_{j}}( b_{j} f) +\rho R( f, x)\right]\mathrm{d} x=0
\end{equation}
Performing integration by parts twice, we have 
\begin{equation} \label{PSgyB5g}
\begin{aligned}
\int _{\Omega }\left[ -\frac{1}{2}\sum _{i,j=1}^{n}\frac{\partial ^{2} \rho }{\partial x_{i} \partial x_{j}} a_{ij} f-\frac{\partial \rho }{\partial x_{j}} b_{j} f+\rho R( f, x)\right]\mathrm{d} x\quad \quad \quad \quad \quad  & \\
+\int _{\partial \Omega }\left[ -\frac{1}{2} \rho \sum _{i,j=1}^{n} \nu _{i}\frac{\partial ( a_{ij} f)}{\partial x_{j}} +\frac{1}{2}\sum _{i,j=1}^{n} \nu _{j}\frac{\partial \rho }{\partial x_{i}} a_{ij} f+\sum _{j=1}^{n} \rho \nu _{j} b_{j} f\right]\mathrm{d} S & =0.
\end{aligned}
\end{equation}
Plugging in the boundary conditions, we get
\begin{equation} \label{aCobixl}
S_{\text{I}} +S_{\text{D}} +S_{\text{N}} =0,
\end{equation}
where
\begin{equation} \label{Fpr20xk}
\begin{cases}
\displaystyle S_{\text{I}} =\int _{\Omega }\left[ -\frac{1}{2}\sum _{i,j=1}^{n}\frac{\partial ^{2} \rho }{\partial x_{i} \partial x_{j}} a_{ij} f-\frac{\partial \rho }{\partial x_{j}} b_{j} f+\rho R( f, x)\right]\mathrm{d} x,\\
\displaystyle S_{\text{D}} =\int _{\partial \Omega _{\text{D}}}\left[\frac{1}{2}\sum _{i,j=1}^{n} \nu _{j}\frac{\partial \rho }{\partial x_{i}} a_{ij} g\right]\mathrm{d} S,\quad S_{\text{N}} =\int _{\partial \Omega _{\text{N}}}\left[ -\rho \phi +\frac{1}{2}\sum _{i,j=1}^{n} \nu _{j}\frac{\partial \rho }{\partial x_{i}} a_{ij} f\right]\mathrm{d} S,
\end{cases}
\end{equation}
The weak formulation given by (\ref{aCobixl}) and (\ref{Fpr20xk}) provides a derivative-free formulation for the function $f$. Note that there are no derivatives of $f$ in the expressions for $S_{\text{I}}$, $S_{\text{D}}$, and $S_{\text{N}}$ in (\ref{Fpr20xk}). Indeed, for a fixed $f$, $S_{\text{I}}$ can be evaluated using a stochastic formulation: if $X_{t}^{x}$ solves the SDE
\begin{equation} \label{tTP3gQa}
\mathrm{d} X_{t}^{x} =b(X_{t}^{x} )\mathrm{d} t+ \sigma (X_{t}^{x} )\mathrm{d} W_{t} ,\quad X_{0}^{x} =x,
\end{equation}
where $W_{t}$ is an $n_{\text{W}}$-dimensional standard Brownian motion, and $\sigma $ is an ($n\times n_{\text{W}}$)-matrix-valued function on $\Omega $ satisfying $\sigma \sigma ^{\top } =a$, then 
\begin{equation} \label{mmG0toJ}
\lim _{\tau \rightarrow 0}\mathbb{E}\frac{\rho (X_{\tau}^{x} )-\rho ( x)}{\tau} =\frac{1}{2}\sum _{i,j=1}^{n}\frac{\partial ^{2} \rho }{\partial x_{i} \partial x_{j}} a_{ij} +\frac{\partial \rho }{\partial x_{j}} b_{j} .
\end{equation}
Indeed, according to Itô's formula, for a fixed $\tau  >0$, we have 
\begin{equation} \label{TC5dit6}
\rho (X_{\tau}^{x} )-\rho ( x) =\sqrt{\tau} \sigma (X_{\tau}^{x} )\xi +\tau \left(\frac{1}{2}\sum _{i,j=1}^{n}\frac{\partial ^{2} \rho }{\partial x_{i} \partial x_{j}} a_{ij} +\frac{\partial \rho }{\partial x_{j}} b_{j} \right) +\mathcal{O}( \tau ^{3/2}) ,
\end{equation}
where $\xi $ is an $n_{\text{W}}$-dimensional standard normal random vector. The directional derivatives of the form $\sum _{i,j=1}^{n} \nu _{j}\frac{\partial \rho }{\partial x_{i}} a_{ij}$ in $S_{\text{D}}$ and $S_{\text{N}}$ are also available through numerical differentiation:
\begin{equation} \label{qbxfRvI}
\sum _{i,j=1}^{n} \nu _{j}\frac{\partial \rho }{\partial x_{i}} a_{ij}\approx \mathcal{D}_{\nu }^{\tilde{\tau}}[ \rho ] (x)=\frac{\rho ( x+\tilde{\tau} A\nu ) -\rho ( x)}{\tilde{\tau}} ,\quad (\tilde{\tau} A\nu )_{j} ={\sum _{j=1}^{n} \nu _{j}( x) a_{ij}}( x) ,
\end{equation}

In this way, we can naturally construct a loss function according to the weak formulation (\ref{aCobixl}). The process is described as follows:
\begin{enumerate}
\item Sample $M_{\text{I}}$ points $(x^{( i)} )_{i=1}^{M_{\text{I}}}$ in the domain $\Omega $, $M_{\text{D}}$ points $(x^{( i)} )_{i=1}^{M_{\text{D}}}$ on the boundary $\partial \Omega _{\text{D}}$, and $M_{\text{N}}$ points $(x^{( i)} )_{i=1}^{M_{\text{N}}}$ on the boundary $\partial \Omega _{\text{N}}$.
\item Compute a numerical approximation of $S_{\text{I}}$: 
\begin{equation} \label{MhBUKFT}
\begin{cases}
\displaystyle \hat{S}_{\text{I}} :=\frac{|\Omega |}{M_{\text{I}}}\sum _{i=1}^{M_{\text{I}}}\left[ -\frac{\rho (X_{\tau}^{x^{( i)}} )-\rho (x^{( i)} )}{\tau} f (x^{(i)})+\rho (x^{( i)} )R\left( f(x^{( i)} )\right)\right] & \\
\displaystyle \hat{S}_{\text{D}} :=\frac{|\partial \Omega _{\text{D}} |}{M_{\text{D}}}\sum _{i=1}^{M_{\text{D}}}\left[\frac{1}{2}\mathcal{D}_{\nu }^{\tilde{\tau}}[ \rho ] (x^{( i)} )g(x^{( i)} )\right] & \\
\displaystyle \hat{S}_{\text{N}} :=\frac{|\partial \Omega _{\text{N}} |}{M_{\text{N}}}\sum _{i=1}^{M_{\text{N}}}\left[ -\rho (x^{( i)} )\phi (x^{( i)} )+\frac{1}{2}\mathcal{D}_{\nu }^{\tilde{\tau}}[ \rho ] (x^{( i)} )f(x^{( i)} )\right] , & 
\end{cases}
\end{equation}

where $X_{\tau}^{x}$ is defined in (\ref{tTP3gQa}) and $\mathcal{D}_{\nu }[ \rho ] (x^{( i)} )$ is defined in (\ref{qbxfRvI}).
\item Construct the loss function
\begin{equation} \label{CCyiKMa}
\mathcal{L} :=\hat{S}_{\text{I}} +\hat{S}_{\text{D}} +\hat{S}_{\text{N}} .
\end{equation}
\end{enumerate}

It should be noted that, in general, before averaging, the difference
\begin{equation*}
\frac{\rho (X_{\tau}^{x} )-\rho (x)}{\tau} \approx \underbrace{\frac{ \sigma (x)\xi }{\sqrt{\tau}}}_{\mathcal O(\tau^{-1/2})} +\underbrace{\left(\frac{1}{2}\sum _{i,j=1}^{n}\frac{\partial ^{2} \rho }{\partial x_{i} \partial x_{j}} a_{ij} +\frac{\partial \rho }{\partial x_{j}} b_{j} \right) (x)}_{\mathcal O(1)}
\end{equation*}
contains a noise term of scale $\tau ^{-1/2}$, which grows to infinity as $\tau \rightarrow 0$. Our method seeks to average out this error by taking the \textit{sample mean over the points} $(x^{( i)} )_{i=1}^{M_{\text{I}}}$. In a rough sense, the error will shrink to $M_{\text{I}}^{-1/2}$ of its original size by averaging $M_{\text{I}}$ samples of noise of scale $\mathcal{O}( 1)$. Therefore, this noise will be of scale $( M_{\text{I}} \tau )^{-1/2}$. To obtain reasonable results, we need 
\begin{equation} \label{1Q9LPhi}
M_{\text{I}}\gg \frac{1}{\tau}
\end{equation}
in order to ``annihilate" the $\tau ^{-1/2}$-scaled error. After constructing the loss, we can naturally apply the weak adversarial training framework by repeating the following two steps in each training epoch:
\begin{enumerate}
\item Fix the parameters of $\rho$, compute $\mathcal{L}$ according to (\ref{CCyiKMa}), and perform gradient descent on the parameters of $f$.
\item Fix the parameters of $f$, compute $\mathcal{L}$ according to (\ref{CCyiKMa}), and perform gradient ascent on the parameters of $\rho $.
\end{enumerate}
In the above framework, we compute the differential operators on $\rho$ using finite difference methods. According to \cite{Chen2025}, replacing differential operators with finite differences usually significantly worsens the results when training deep neural networks to solve PDEs. Nevertheless, in our framework, $\rho$ is the test function—we do not really care about its accuracy (indeed, we cannot even define ``accuracy", since there is no such thing as a ``true test function").  

\subsection{Weakly Nonlinear Parabolic PDEs}
A general parabolic Fokker-Planck type PDE involving time reads
\begin{equation}
    \label{teFP}
    \frac{\partial f}{\partial t}-\frac{1}{2}\sum _{i,j=1}^{n}\frac{\partial ^{2}( a_{ij} f)}{\partial x_{i} \partial x_{j}} +\sum _{j=1}^{n}\frac{\partial }{\partial x_{j}}( b_{j} f) +R( f, x, t) =0,\quad (x,t)\in \Omega \times (0,T),
\end{equation}
with initial condition
\begin{equation}
    f(x, 0) = f_0(x),\quad x \in \overline \Omega
\end{equation}
and boundary conditions 
\begin{equation}
\begin{cases}
f=g & (x,t)\in (0, T)\times \partial \Omega _{\text{D}}\\
\displaystyle \sum _{j=1}^{n} \nu _{j}\left[\frac{1}{2}\sum _{i=1}^{n}\frac{\partial }{\partial x_{i}} (a_{ij} f)-b_{j} f\right] =\phi  & (x,t)\in (0, T)\times \partial \Omega _{\text{N}}
\end{cases}
\end{equation}
For this type of equation, the weak form is very similar to the elliptic case, and thus the method is nearly identical.

\subsection{Second-Order PDEs not of Fokker-Planck Type}
If the second-order operator takes other forms, e.g.,
\begin{equation*}
\frac{1}{2} a_{ij}\sum _{i,j=1}^{n}\frac{\partial ^{2} f}{\partial x_{i} \partial x_{j}} +\sum _{j=1}^{n} b_{j}\frac{\partial f}{\partial x_{j}} +R( f, x) =0
\end{equation*}
or 
\begin{equation*}
-\frac{1}{2}\sum _{i,j=1}^{n}\frac{\partial }{\partial x_{i}}\left( a_{ij}\frac{\partial f}{\partial x_{j}}\right) +\sum _{j=1}^{n} b_{j}\frac{\partial f}{\partial x_{j}} +R( f, x) =0,
\end{equation*}
the framework will still work. This will involve computing the derivatives of the products of the coefficients and the test functions, i.e., $a_{ij} \rho$ and $b_j \rho$. Still, the formulation remains derivative-free for $f$. The same applies to their parabolic counterparts:
\begin{equation*}
\frac{\partial f}{\partial t} + \frac{1}{2} a_{ij}\sum _{i,j=1}^{n}\frac{\partial ^{2} f}{\partial x_{i} \partial x_{j}} +\sum _{j=1}^{n} b_{j}\frac{\partial f}{\partial x_{j}} +R( f, x, t) =0
\end{equation*}
or 
\begin{equation*}
\frac{\partial f}{\partial t} - \frac{1}{2}\sum _{i,j=1}^{n}\frac{\partial }{\partial x_{i}}\left( a_{ij}\frac{\partial f}{\partial x_{j}}\right) +\sum _{j=1}^{n} b_{j}\frac{\partial f}{\partial x_{j}} +R( f, x, t) =0.
\end{equation*}

\section{Conclusion}
In this paper, we proposed ARDO, which provides an alternative to DRDM, DeepMartNet and other stochastic neural PDE solvers. This method is more robust since it is derivative-free with respect to the solution network. Future numerical experiments will be conducted to evaluate the efficiency and accuracy of the proposed algorithm, especially in high dimensions.

\end{document}